\documentclass[12pt,a4paper]{article}

\usepackage{amsfonts, amsmath, amsthm, amssymb, bm, dsfont}
\usepackage{algorithm, algorithmic}  
\usepackage{stmaryrd}  
\usepackage{graphicx}
\usepackage{path}
\usepackage{datetime}  

\usepackage[paper=a4paper,dvips,top=2cm,left=2cm,right=2cm,foot=1cm,bottom=3cm]{geometry}

\begin{document}

\title{\Large\bf A Structured Levenberg--Marquardt Method for \\ Binary Tensor Completion}

\author{
  Yannan Chen\footnote{School of Mathematical Sciences, and Key Laboratory of Mathematical Foundation and Computation in Nonlinear Sciences of Guangdong Higher Education Institutes, South China Normal University, Guangzhou 510631, Guangdong Province, China. (\texttt{email: ynchen@scnu.edu.cn})} \and
  Xihong Yan\footnote{Shanxi Key Laboratory for Intelligent Optimization Computing and Block-chain Technology, School of Mathematics and Statistics, Taiyuan Normal University, Jinzhong 030619, Shanxi Province, China (\texttt{email: yanxihong@tynu.edu.cn}).} \and
  Bobin Zheng\footnote{School of Mathematical Sciences, South China Normal University, Guangzhou 510631, Guangdong Province, China (\texttt{email: 2022021959@m.scnu.edu.cn}).}
 }

\date{\today}
\maketitle

\begin{abstract}
We develop a structured Levenberg--Marquardt (LM) method with trust-region globalization for binary tensor completion, using triple decomposition and a logistic loss. Under a local component-residual error bound and positive outer curvature, we prove quadratic convergence of the exact iteration for damping exponents $1\le\delta\le2$, without requiring full row rank. A constant-rank criterion and an explicit rank-one certificate provide sufficient conditions for the residual error bound. The analysis decomposes the local error into outer-loss linearization, factor linearization, and regularization contributions. For a symmetric nonquadratic logistic--LSE family with $\delta=2$, cancellation improves the objective-model error to fourth order and the decrease-ratio error to second order. A synthetic experiment illustrates the local quadratic convergence rate. On the HCP and Enron datasets, the proposed method improves binary tensor completion accuracy while requiring fewer outer iterations than the first-order methods.

		
\end{abstract}

	\newtheorem{theorem}{Theorem}[section]
	\newtheorem{lemma}[theorem]{Lemma}
	\newtheorem{definition}[theorem]{Definition}
	\newtheorem{assumption}[theorem]{Assumption}
	\newtheorem{example}[theorem]{Example}
	\newtheorem{remark}[theorem]{Remark}
\newtheorem{proposition}[theorem]{Proposition}
\newtheorem{corollary}[theorem]{Corollary}

	\newcommand{\REAL}{\mathds{R}}
	\newcommand{\Tt}{\Theta}
	\newcommand{\At}{\mathcal{A}}
	\newcommand{\Bt}{\mathcal{B}}
	\newcommand{\Ct}{\mathcal{C}}
	\newcommand{\Am}{\mathbf{A}}
	\newcommand{\Bm}{\mathbf{B}}
	\newcommand{\Cm}{\mathbf{C}}
	\newcommand{\Dm}{\mathbf{D}}
	\newcommand{\Em}{\mathbf{E}}
	\newcommand{\Om}{\mathbf{O}}
	\newcommand{\Um}{\mathbf{U}}
	\newcommand{\Vm}{\mathbf{V}}
	\newcommand{\Mm}{\mathbf{M}}
	\newcommand{\Pm}{\mathbf{P}}
	\newcommand{\Xm}{\mathbf{X}}
	\newcommand{\Zm}{\mathbf{Z}}
	\newcommand{\Ten}{\mathcal{T}}
	\newcommand{\Yt}{\mathcal{Y}}
	\newcommand{\Zt}{\mathcal{Z}}
	\newcommand{\av}{\mathbf{a}}
	\newcommand{\bv}{\mathbf{b}}
	\newcommand{\cv}{\mathbf{c}}
	\newcommand{\gv}{\mathbf{g}}
	\newcommand{\pv}{\mathbf{p}}
	\newcommand{\sv}{\mathbf{s}}
	\newcommand{\rv}{\mathbf{r}}
	\newcommand{\yv}{\mathbf{y}}
	\newcommand{\xv}{\mathbf{x}}
	\newcommand{\zv}{\mathbf{z}}
	\newcommand{\thetav}{{\bm\theta}}

	\newcommand{\dist}{\mathrm{dist}}
	\renewcommand{\vec}{\mathrm{vec}}
	\newcommand{\Sigmam}{\mathbf{\Sigma}}
	\newcommand{\Jac}{\mathbf{J}}
	\newcommand{\Hess}{\mathbf{H}}
	\newcommand{\Pcd}{\mathbf{M}}
	\newcommand{\norm}[1]{\left\|#1\right\|}
	\newcommand{\diff}{\mathrm{d}}
	\newcommand{\Z}{\mathds{Z}}
	\newcommand{\ten}[1]{\mathcal{#1}}
	\newcommand{\Xt}{\mathcal{X}}
	\newcommand{\g}{\boldsymbol{g}}
	\newcommand{\p}{\boldsymbol{p}}
	\newcommand{\z}{\boldsymbol{z}}
	\newcommand{\reflack}{{ (reference)}}
	\newcommand{\argmin}{\mathrm{argmin}}
	\newcommand{\Diag}[1]{\mathrm{diag}\left(#1\right)}
	\newcommand{\blkDiag}[1]{\mathrm{blkdiag}\left(#1\right)}
	\newcommand{\bZ}{{\bf Z}}
	\newcommand{\ov}[1]{\overline{#1}}

	\newpage

\section{Introduction}\label{Sec-1}

Binary tensors describe multiway observations of whether an event occurs, including brain connectivity across subjects and communication between users over time. Completion predicts missing outcomes through shared low-rank latent factors. This representation preserves interactions among the modes and allows observed entries to inform predictions at missing locations. Because observations take only two values, a probability model distinguishes latent signal strength from realized outcomes. Probabilistic factorization and likelihood-based binary tensor models provide a natural framework \cite{RHHC15,WL20,ZL21}. Rai et al.\ \cite{RHHC15} develop scalable probabilistic factorization for binary and count data. Wang and Li \cite{WL20} establish a multilinear Bernoulli model and statistical accuracy guarantees, clarifying how low-rank structure and binary observations affect estimation. These developments motivate algorithms that exploit both the curvature of the binary loss and the structure of the low-rank tensor representation.

The tensor representation determines both the model class and the optimization variables. Generalized canonical polyadic (CP) fitting \cite{HKD20} accommodates non-least-squares losses, allowing the objective to match the data distribution. The triple decomposition of Qi et al.\ \cite{QCBZ21} represents a third-order tensor by contracting three factor tensors over three internal indices. Each factor participates in two internal modes, giving explicit multilinear derivatives that can be exploited in optimization. Its structure has motivated several optimization methods: Cui and Chang \cite{CC24} apply limited-memory BFGS (L-BFGS) to hyperspectral super-resolution, Chen et al.\ \cite{CLLX25} develop composite generalized Gauss--Newton (GGN) for Poisson completion, and Wen et al.\ \cite{WLCC26} study Riemannian preconditioning for low-rank completion. These works demonstrate the value of combining the triple representation with suitable losses and optimization methods.

We model binary observations using a logistic likelihood with a log-sum-exp (LSE) margin penalty. The likelihood promotes agreement with observed labels, while the penalty controls excessively large positive signed margins. The resulting objective is a separable loss composed with a trilinear tensor map. Its diagonal outer Hessian and explicit factor derivatives are well suited to structured second-order optimization. Vandecappelle et al.\ \cite{VVL20} develop GGN for non-least-squares CP fitting and analyze local quadratic behavior. Their inexact extension \cite{VVL21} uses implicit derivative products to improve scalability. Following this approach, we derive forward and adjoint products at the observed entries of the triple model and construct a factor-row block preconditioner. The implementation avoids forming the full Jacobian; operation counts describe the costs of contractions, preconditioning, and inner iterations. 

Local error-bound theory supports fast convergence analysis for factorized models with nonisolated solutions. Yamashita and Fukushima \cite{YF01} establish quadratic Levenberg--Marquardt (LM) convergence under a local residual error bound. Fan and Yuan \cite{FY05} allow damping exponents in $[1,2]$ and prove quadratic convergence of the iterates without Jacobian nonsingularity. Bergou et al.\ \cite{BDK20} develop further convergence and complexity results, while Fischer et al.\ \cite{FIS24} review LM methods for singular problems. For weighted composite models, Wu et al.\ \cite{WHZW24} connect GGN with regularized natural gradient methods and establish local quadratic convergence using a component-loss residual and an error bound. The error-bound viewpoint is particularly useful when scaling freedoms make a solution nonunique, as often occurs in tensor factorization. These results provide tools for treating the redundant factorization and variable curvature in binary triple completion.

Here, $G$ collects the outer-loss derivatives at fitted observations, $J$ is the triple-map Jacobian, and $H$ is the diagonal outer Hessian. The exact LM system uses $J^THJ+\|G\|^\delta E$, $\delta\in[1,2]$. Since $\mathrm{D}G=HJ$, where $\mathrm{D}$ denotes differentiation, the analysis accounts for variations in the outer Hessian and for singular directions of the tensor-map Jacobian. A weighted Taylor comparison and a separation of persistent and small singular directions establish quadratic convergence near a realizable root under a component-residual error bound and positive outer curvature. Growth conditions provide a complementary interpretation. Marumo et al.\ \cite{MOT25} develop a generalized LM method with oracle complexity guarantees and local quadratic convergence under suitable growth and zero-gap conditions. Rebjock and Boumal \cite{RB25} establish relationships among quadratic growth, gradient error bounds, the Polyak--\L{}ojasiewicz condition, and the Morse--Bott property for smooth objectives. In our realizable positive-curvature setting, the component-residual error bound is equivalent to local quadratic growth, connecting the local assumption to objective geometry and the distance to the solution set.

Rank-one completion provides concrete verification tools. Kahle et al.\ \cite{KKKR17} study its algebraic and geometric structure, including exponent matrices and real-completion sign constraints. Zhou et al.\ \cite{ZNPZ25} study cubic rank-one completion; subsequent work develops stable noisy completion through biquadratic optimization \cite{NTZ26} and recursive spectral methods \cite{ZDN26}. Using the sign and exponent-matrix description, we obtain an explicit error-bound certificate for realizable rank-one binary models. Its constant links loss curvature, target magnitude, observation-pattern singular values, and factor scale. The certificate therefore makes the local assumption interpretable in terms of both the loss and the sampling structure.

The logistic--LSE loss also permits a refinement of model agreement. For a symmetric parameter family, its third derivative vanishes at the scalar minimizer. Combined with an exact LM step with $\delta=2$, this improves the objective-model discrepancy from cubic to fourth order and the decrease-ratio error from first to second order, while preserving quadratic iterate convergence.

The main contributions are summarized as follows.
\begin{enumerate}
	\item We develop a structured LM algorithm with trust-region globalization, observed-entry forward and adjoint contractions, and factor-row block-diagonal preconditioning. We give operation counts and verify convergence to stationarity for bounded iterates and sufficient model decrease.
	\item We establish local quadratic convergence of exact original-weight steps for every $\delta\in[1,2]$ under a component-residual error bound and positive outer curvature. The analysis allows rank-deficient and rank-changing Jacobians, and shows that the trust-region method accepts exact full steps locally with the original curvature. A constant-rank criterion and an explicit rank-one certificate provide verification routes.
	\item We separate outer-loss, mixed-factor, and damping errors, and obtain fourth-order model error and second-order decrease-ratio agreement for a symmetric nonquadratic logistic--LSE family. The synthetic experiment illustrates quadratic local convergence. Experiments on HCP and Enron compare prediction accuracy, outer iteration counts, and computational time.
\end{enumerate}

The remainder of the paper is organized as follows. Section~\ref{Sec-2} presents the model, algorithm, and computational implementation. Section~\ref{Sec-3} introduces the residual error bound and verification criteria. Section~\ref{Sec-4} establishes the error decomposition, local quadratic convergence, and model-agreement refinement. Section~\ref{Sec-7} presents the synthetic and real-data experiments, and Section~\ref{Sec-8} concludes the paper.

\section{Binary triple model and structured LM algorithm}\label{Sec-2}

\subsection{Penalized logistic model}

Let $I,J,K$ be the three mode sizes. The unknown real tensor $\Theta=(\theta_{ijk})$ is a deterministic parameter: its entries describe latent signal strengths, rather than binary outcomes or probabilities. The binary response $Y_{ijk}\in\{0,1\}$ is a random variable, and $y_{ijk}$ denotes its realized observed value. We fit the model using the index set $\Omega\subset\{1,\ldots,I\}\times\{1,\ldots,J\}\times\{1,\ldots,K\}$, with $n=|\Omega|$. An observed value $y_{ijk}=0$ is a label; an index outside $\Omega$ is missing. All probabilities below refer to the model at a fixed value of $\Theta$.

Following the latent-threshold formulation in \cite{WL20}, we model each response by adding random noise to the latent entry and recording whether the result is nonnegative. Specifically, the noise variables $Z_{ijk}$ are independent and have a logistic distribution with location parameter $0$ and fixed scale parameter $\beta>0$. This distribution is specified by its cumulative distribution function
\[
 p(u):=\mathbb P(Z_{ijk}\le u)
      =\frac{1}{1+e^{-u/\beta}},\qquad u\in\REAL.
\]
Thus $p(u)$ is the probability that the noise is at most $u$. The distribution is symmetric about $0$, has mean $0$, and satisfies $1-p(-u)=p(u)$. Given $\Theta$, the observation model is
\begin{equation}\label{model:latVar}
 Y_{ijk}=\mathbf1\{\theta_{ijk}+Z_{ijk}\ge0\},\qquad
 \mathbb P(Y_{ijk}=1\mid\Theta)=p(\theta_{ijk}).
\end{equation}
The indicator equals $1$ when signal plus noise is nonnegative and $0$ otherwise. To obtain the second relation, hold $\theta_{ijk}$ fixed. Since the noise distribution is continuous,
\[
 \mathbb P(Y_{ijk}=1\mid\Theta)
 =\mathbb P(Z_{ijk}\ge-\theta_{ijk})
 =1-p(-\theta_{ijk})=p(\theta_{ijk}).
\]
Consequently, $Y_{ijk}$ has a Bernoulli distribution: its probabilities of being $1$ and $0$ are $p(\theta_{ijk})$ and $1-p(\theta_{ijk})$, respectively. In particular, a positive latent entry makes label $1$ more likely, a negative entry makes label $0$ more likely, and a zero entry gives equal probabilities. For fixed nonzero $\theta_{ijk}$, increasing $\beta$ moves the probability toward $1/2$. We fix $\beta$ during fitting because the probabilities depend only on $\theta_{ijk}/\beta$; multiplying both $\Theta$ and $\beta$ by the same positive constant leaves the observation model unchanged.

Independence of the noise variables implies independence of the responses for a fixed $\Theta$. The probability of the observed labels is therefore a product of Bernoulli probabilities. Taking its logarithm gives the log-likelihood, viewed as a function of $\Theta$ with the observed labels held fixed:
\[
 \ell(\Theta;\mathbf y_\Omega)=\sum_{(i,j,k)\in\Omega}
 \{y_{ijk}\log p(\theta_{ijk})+(1-y_{ijk})\log[1-p(\theta_{ijk})]\}.
\]
For label $1$, the summand is $\log p(\theta_{ijk})$; for label $0$, it is $\log[1-p(\theta_{ijk})]$. Define the observed signed label $s_{ijk}:=2y_{ijk}-1$, which equals $1$ for label $1$ and $-1$ for label $0$. The negative log-likelihood then has the single expression
\[
 -\ell(\Theta;\mathbf y_\Omega)
 =\sum_{(i,j,k)\in\Omega}\log(1+e^{-s_{ijk}\theta_{ijk}/\beta}).
\]
Minimizing this loss is equivalent to maximizing the likelihood. The signed margin $s_{ijk}\theta_{ijk}$ is positive when the sign of the latent entry agrees with the observed label. This likelihood-based formulation is standard in binary tensor fitting \cite{WL20,ZL21}.

We represent the latent tensor by the triple decomposition of \cite{QCBZ21}:
\begin{equation}\label{eq:TTD}
 \theta_{ijk}={\sum_{r_1=1}^R\sum_{r_2=1}^R\sum_{r_3=1}^R}
 a_{ir_2r_3}b_{r_1jr_3}c_{r_1r_2k},
 \qquad \Theta=\llbracket\At,\Bt,\Ct\rrbracket,
\end{equation}
where $\At\in\REAL^{I\times R\times R}$, $\Bt\in\REAL^{R\times J\times R}$, and $\Ct\in\REAL^{R\times R\times K}$ are deterministic optimization variables. The indices $r_1,r_2,r_3$ are summed internal indices, and $R$ is the chosen representation size. The smallest such $R$ for a given tensor is its triple rank. The double-bracket notation denotes exactly the contraction in \eqref{eq:TTD}.

To discourage excessively large fitted margins, we add a smooth positive-part penalty. The standard log-sum-exp function \cite[Section 3.1.5]{BV04} gives
\[
 (t-\alpha)_+\le \frac1\tau\log(1+e^{\tau(t-\alpha)})
 \le(t-\alpha)_++\frac{\log2}{\tau},
 \qquad (u)_+:=\max\{u,0\}.
\]
Here $\alpha>0$ is the margin threshold and $\tau>0$ controls the sharpness of the smooth approximation. The scalar logistic--LSE loss is defined by
\begin{equation}\label{eq:q-calibration}
 q(t):=\log(1+e^{-t/\beta})
      +\frac\lambda\tau\log(1+e^{\tau(t-\alpha)}).
\end{equation}
The first term in \eqref{eq:q-calibration} is the negative Bernoulli log-likelihood for a signed margin $t$. The second term is a smooth margin penalty with weight $\lambda>0$; it penalizes large positive margins and ensures a finite scalar minimizer.

We fit the triple factors by solving
\begin{equation}\label{ourModel}
	\min_{\At,\Bt,\Ct}
	\sum_{(i,j,k)\in\Omega}q(s_{ijk}\theta_{ijk}),
\end{equation}
where $s_{ijk}=2y_{ijk}-1\in\{-1,1\}$ is the fixed signed label. The real number $\theta_{ijk}$ is the latent entry generated by the factors $\At,\Bt,\Ct$ through the triple contraction \eqref{eq:TTD}; it varies with the optimization variables and is not a binary label or a probability. Their product $s_{ijk}\theta_{ijk}$ is the signed margin. Model \eqref{ourModel} sums these losses over the observed set $\Omega$. The chosen size $R$ restricts the representation and couples observed and missing entries through shared factors. 

Differentiating the scalar loss gives
\[
 q'(t)=-\frac{1}{\beta(1+e^{t/{\beta}})}
       +\frac{\lambda}{1+e^{-\tau(t-\alpha)}},
\]
\[
 q''(t)=\frac{e^{t/{\beta}}}{{\beta}^2(1+e^{t/{\beta}})^2}
       +\frac{\lambda\tau e^{\tau(t-\alpha)}}{(1+e^{\tau(t-\alpha)})^2}>0.
\]
The inequality $u/(1+u)^2\le1/4$ for $u>0$ gives
\[
 0<q''(t)\le\frac1{4\beta^2}+\frac{\lambda\tau}{4},
 \qquad \lim_{|t|\to\infty}q''(t)=0.
\]
Thus $q$ is analytic and strictly convex, and a uniform positive lower curvature bound is available on each compact interval. Moreover,
\[
 \lim_{t\to-\infty}q'(t)=-1/\beta<0,
 \qquad \lim_{t\to+\infty}q'(t)=\lambda>0.
\]
By continuity, $q'$ vanishes at a finite point $t_0$; since $q''>0$, this point is unique. The derivative is negative before $t_0$ and positive after it, so $t_0$ is the unique global minimizer of $q$. The local analysis uses the compact-interval curvature bounds stated above. For classification, we choose $q'(0)<0$, equivalently $\lambda<(1+e^{\tau\alpha})/(2\beta)$, so that the preferred signed margin $t_0$ is positive. All fitted parameter choices below satisfy this condition.

\subsection{Composite derivatives and damped step}
All vector norms are Euclidean and all matrix norms are induced operator norms, unless a weight is stated. The symbols $E$ and $E_r$ denote identity matrices of the required size and of size $r$, respectively. Throughout, $\vec$ stacks matrix and tensor entries with the leftmost index varying fastest. Set
\[
 \xv=((\vec\At)^T,(\vec\Bt)^T,(\vec\Ct)^T)^T\in\REAL^d,
 \quad d=(I+J+K)R^2,\quad
 \thetav(\xv)=\Pi_\Omega\vec\llbracket\At,\Bt,\Ct\rrbracket.
\]
The vector $\xv$ contains all factor entries, and $\thetav(\xv)\in\REAL^n$ contains the fitted latent entries at observed indices. Fix an ordering $\Omega=\{\omega_1,\ldots,\omega_n\}$. The coordinate-selection matrix $\Pi_\Omega\in\REAL^{n\times IJK}$ extracts those entries in that order. Thus $\theta_\nu(\xv)=\theta_{\omega_\nu}(\xv)$ and $s_\nu=2y_{\omega_\nu}-1$. Define the univariate component loss $f_\nu(u):=q(s_\nu u)$, so $f_\nu'(u)=s_\nu q'(s_\nu u)$ and $f_\nu''(u)=q''(s_\nu u)$. The separable composite objective and its component-gradient residual are
\begin{equation}\label{ComOptim}
 F(\xv):=\sum_{\nu=1}^nf_\nu(\theta_\nu(\xv)),\qquad
 G(\xv):=(f_\nu'(\theta_\nu(\xv)))_{\nu=1}^n{.}
\end{equation}
Here $G(\xv)\in\REAL^n$ measures how far the observed latent entries are from minimizing their individual losses. It differs from the factor gradient $\nabla F\in\REAL^d$. 
In particular,
\begin{equation}\label{eq:realizable-target}
 G(\xv)=0\quad\Longleftrightarrow\quad
 \theta_\nu(\xv)=s_\nu t_0\quad(\nu=1,\ldots,n).
\end{equation}
This is an interpolation condition on the realized labels.

Define
\[
 J:=\Jac(\xv)=\mathrm{D}\thetav(\xv)\in\REAL^{n\times d},\qquad
 H:=H(\xv)=\Diag{f_\nu''(\theta_\nu(\xv))}\in\REAL^{n\times n}.
\]
Thus $J$ describes the inner tensor map, whereas $G$ and $H$ are the gradient and Hessian of the outer separable loss evaluated at $\thetav(\xv)$. The diagonal positive matrix $H$ contains the curvature of the outer losses; it is not the full Hessian with respect to the factors. The chain rule gives the standard composite derivative formulas \cite{HKD20,VVL20}:
\begin{equation}\label{eq:exact-hessian}
 \nabla F=J^TG,\qquad
 \nabla^2F=J^THJ+\sum_{\nu=1}^nG_\nu\nabla^2\theta_\nu.
\end{equation}
The first term $J^THJ$ is the positive semidefinite GGN approximation. The second term contains the second derivatives of the tensor map and can be indefinite. 

We retain the first positive semidefinite term and add LM regularization \cite{FY05,FIS24}. With subscript $t$ denoting evaluation at iterate $\xv_t$, the damped step solves
\begin{equation}\label{eq:LM method_system}
 (J_t^TH_tJ_t+\mu_tE)\sv_t=-J_t^TG_t,
 \qquad \mu_t=\|G_t\|^\delta,\quad1\le\delta\le2.
\end{equation}
Here $\mu_t$ is the LM regularization parameter and $\delta$ is the regularization exponent. Since $\mu_t>0$ when $G_t\ne0$, adding $\mu_t E$ makes the coefficient matrix positive definite even when $J_t$ is rank deficient. The exact local iteration is
\begin{equation}\label{eq:LM method}
 \xv_{t+1}=\xv_t+\sv_t;
\end{equation}
The algorithm terminates when $\nabla F(\xv_t)=0$, which includes $G_t=0$. The matrix $H_t$ represents the current outer-loss curvature, and $\mu_tE$ regularizes the factor-space system. Theorem~\ref{thm:Quad rate} analyzes exact full steps. The trust-region method below provides globalization.

\subsection{Trust-region globalization}\label{subsec:trust-implementation}\label{Sec-5}
Set $\gv_t:=\nabla F(\xv_t)=J_t^TG_t$. At a nonstationary iterate, a trial step is obtained by approximately solving
\begin{equation}\label{eq:LM method_trust region}
 \min_{\|\sv\|\le\Delta_t}\ \phi_t(\sv)
 :=\gv_t^T\sv+\tfrac12\sv^TB_t\sv,
 \qquad B_t=J_t^TH_t^{\rm tr}J_t+\mu_t E,
 \qquad \mu_t=\|G_t\|^\delta.
\end{equation}
Here $\Delta_t\in(0,\Delta_{\max}]$ is the trust-region radius. The original-curvature model uses $H_t^{\rm tr}=H_t$. In the real-data implementation, small diagonal entries of $H_t$ are replaced by a fixed positive floor: $(H_t^{\rm tr})_{\nu\nu}=\max\{(H_t)_{\nu\nu},h_{\rm floor}\}$. This changes only the quadratic model; both $F$ and $\gv_t$ are evaluated as before. Since $\gv_t\ne0$ implies $G_t\ne0$, the matrix $B_t$ is positive definite in either case.

We use preconditioned Steihaug CG \cite{S83,NW06}, which truncates an inner step at the Euclidean trust-region boundary. Write $P_t=-\phi_t(\sv_t)$ for the predicted reduction and require
\begin{equation}\label{eq:cauchy-decrease}
 P_t\ge \kappa_{\rm cd}\|\gv_t\|
    \min\left\{\Delta_t,\frac{\|\gv_t\|}{1+\|B_t\|}\right\},
 \qquad 0<\kappa_{\rm cd}\le\tfrac12 .
\end{equation}
For example, one may require at least half the Cauchy reduction and use the Cauchy step if this test fails. The Cauchy step minimizes $\phi_t$ on the segment in the direction $-\gv_t/\|\gv_t\|$ of length $\Delta_t$; at the shorter of $\Delta_t$ and $\|\gv_t\|/(1+\|B_t\|)$, its model reduction is at least half the product of that length and $\|\gv_t\|$. Thus the fallback enforces \eqref{eq:cauchy-decrease} with $\kappa_{\rm cd}=1/4$. In particular, $P_t>0$, and the acceptance ratio is well defined:
\begin{equation}\label{rho-def}
 \rho_t=\frac{F(\xv_t)-F(\xv_t+\sv_t)}{P_t}.
\end{equation}

\begin{algorithm}[t]
\caption{Structured trust-region LM for binary triple completion}\label{alg:TRCGLM}
\begin{algorithmic}[1]
\REQUIRE $\xv_0$, $0<\Delta_0\le\Delta_{\max}$, $\eta_{\rm a}\in(0,1/4)$, $\delta\in[1,2]$, tolerance $\epsilon_{\rm opt}\ge0$, and limit $T_{\max}$.
\STATE Set $t\gets0$ and evaluate $\gv_0$.
\WHILE{$\|\gv_t\|>\epsilon_{\rm opt}$ and $t<T_{\max}$}
\STATE Form the model \eqref{eq:LM method_trust region} using structured products.
\STATE Compute a feasible trial step by truncated CG, enforcing \eqref{eq:cauchy-decrease} by the Cauchy fallback if needed.
\STATE Evaluate $\rho_t$ by \eqref{rho-def}.
\IF{$\rho_t\ge\eta_{\rm a}$}
  \STATE Set $\xv_{t+1}=\xv_t+\sv_t$.
\ELSE
  \STATE Set $\xv_{t+1}=\xv_t$.
\ENDIF
\IF{$\rho_t<1/4$}
  \STATE Set $\Delta_{t+1}=\Delta_t/4$.
\ELSIF{$\rho_t>3/4$ and $\|\sv_t\|=\Delta_t$}
  \STATE Set $\Delta_{t+1}=\min\{2\Delta_t,\Delta_{\max}\}$.
\ELSE
  \STATE Set $\Delta_{t+1}=\Delta_t$.
\ENDIF
\STATE Set $t\gets t+1$ and evaluate $\gv_t$.
\ENDWHILE
\RETURN $\xv_t$.
\end{algorithmic}
\end{algorithm}

\begin{theorem}[Global convergence {\cite[Theorem~4.6]{NW06}}]\label{prop:TR-global}
Suppose that the iterates of Algorithm~\ref{alg:TRCGLM} are bounded and that the feasible trial steps satisfy \eqref{eq:cauchy-decrease} with fixed $\kappa_{\rm cd}>0$. With $\epsilon_{\rm opt}=0$ and $T_{\max}=\infty$, the algorithm either terminates at a stationary point or generates a sequence for which $\|\nabla F(\xv_t)\|\to0$. Every accumulation point is stationary. 
\end{theorem}

This is the classical first-order convergence result for trust-region methods with a positive acceptance threshold; see \cite[Theorem~4.6]{NW06}. In our setting, $F$ is bounded below, and bounded iterates and the finite radius cap give a uniform Lipschitz bound for $\nabla F$ on the trial segments and a uniform bound on $\|B_t\|$. Thus \eqref{eq:cauchy-decrease} supplies the sufficient model reduction used in the standard convergence argument.


\subsection{Observed-entry contractions and preconditioning}
The following formulas specialize the derivative-product principle of \cite{VVL20,VVL21} to the triple contraction. At a fixed iterate, suppress the iteration subscript and write $B=J^TH^{\rm tr}J+\mu E$. For a factor perturbation $\mathbf d=((\vec\Delta\At)^T,(\vec\Delta\Bt)^T,(\vec\Delta\Ct)^T)^T$, trilinearity gives the forward product
\begin{equation}\label{eq:HJx}
\begin{split}
 (J\mathbf d)_{ijk}={\sum_{r_1=1}^R\sum_{r_2=1}^R\sum_{r_3=1}^R}\big(&
 \Delta a_{ir_2r_3}b_{r_1jr_3}c_{r_1r_2k}
 +a_{ir_2r_3}\Delta b_{r_1jr_3}c_{r_1r_2k}\\
 &+a_{ir_2r_3}b_{r_1jr_3}\Delta c_{r_1r_2k}\big),
 \quad(i,j,k)\in\Omega.
\end{split}
\end{equation}
For example, differentiating $\theta_{ijk}$ with respect to $a_{ir_2r_3}$ leaves a contraction of $\Bt$ and $\Ct$. 

Define these pairwise derivative coefficients by
\[
 d^A_{jk;r_2r_3}:={\sum_{r_1=1}^R}b_{r_1jr_3}c_{r_1r_2k},\quad
 d^B_{ik;r_1r_3}:={\sum_{r_2=1}^R}a_{ir_2r_3}c_{r_1r_2k},\quad
 d^C_{ij;r_1r_2}:={\sum_{r_3=1}^R}a_{ir_2r_3}b_{r_1jr_3}.
\]
For $\mathbf z\in\REAL^n$, let $z_{ijk}$ denote its component associated with the observed index $(i,j,k)$. The adjoint product $J^T\mathbf z$ accumulates weighted derivatives into arrays $\ten Z_A,\ten Z_B,\ten Z_C$ of the same sizes as the factors:
\begin{equation}\label{eq:JTz}
\begin{split}
 ({\ten Z_A})_{ir_2r_3}&:={\sum_{j,k:\,(i,j,k)\in\Omega}}z_{ijk}d^A_{jk;r_2r_3},\\
 ({\ten Z_B})_{r_1jr_3}&:={\sum_{i,k:\,(i,j,k)\in\Omega}}z_{ijk}d^B_{ik;r_1r_3},\\
 ({\ten Z_C})_{r_1r_2k}&:={\sum_{i,j:\,(i,j,k)\in\Omega}}z_{ijk}d^C_{ij;r_1r_2}.
\end{split}
\end{equation}
The physical index on the left is fixed in each sum, and $J^T\mathbf z=((\vec\ten Z_A)^T,(\vec\ten Z_B)^T,(\vec\ten Z_C)^T)^T$. Thus $B\mathbf d$ is evaluated through $J\mathbf d$, diagonal weighting, and $J^T$ accumulation, followed by $\mu\mathbf d$. No dense tensor of size $IJK$ is required.

Partition the Jacobian as $J=[J_A\;J_B\;J_C]$, where each block differentiates with respect to one entire factor tensor. A block-Jacobi preconditioner retains the diagonal factor blocks of the regularized GGN matrix and omits cross-factor couplings. Following \cite{VVL20,VVL21}, use
\[
 \Pm=\operatorname{diag}(J_A^TH^{\rm tr}J_A,
                         J_B^TH^{\rm tr}J_B,
                         J_C^TH^{\rm tr}J_C)+\mu E.
\]
After a fixed coordinate permutation, each factor block is block diagonal with respect to the physical index. For example,
\[
 \Pm_{A,i}:={\sum_{j,k:\,(i,j,k)\in\Omega}}h^{\rm tr}_{ijk}
       \mathbf d^A_{jk}(\mathbf d^A_{jk})^T+\mu {E_{R^2}},
\]
where $h^{\rm tr}_{ijk}$ is the corresponding diagonal entry of $H^{\rm tr}$, and $\mathbf d^A_{jk}$ stacks $d^A_{jk;r_2r_3}$ with $r_2$ changing fastest. The blocks $\Pm_{B,j}$ and $\Pm_{C,k}$ are obtained by summing over $i,k:(i,j,k)\in\Omega$ and $i,j:(i,j,k)\in\Omega$, respectively. 

A direct forward/adjoint pair costs $\mathcal O(nR^3+d)$ with $\mathcal O(n+d)$ storage. {If $N_{\rm pair}$ is the total number of distinct observed pairs across the three modes,} caching pair contractions costs ${\mathcal O(N_{\rm pair}R^3)}$ work and ${\mathcal O(N_{\rm pair}R^2)}$ storage; each subsequent GGN product costs $\mathcal O(nR^2+d)$. With $N_{\rm CG}$ products, the contraction work is therefore
\[
 \mathcal O\bigl({N_{\rm pair}}R^3+N_{\rm CG}(nR^2+d)\bigr).
\]
Straightforward preconditioner setup costs $\mathcal O(nR^4)$, its block factorizations cost $\mathcal O((I+J+K)R^6)$, and an application costs $\mathcal O((I+J+K)R^4)$. 

\section{Local residual error bounds}\label{Sec-3}

\begin{definition}[Local residual error bound]\label{LEBnew}
Let $\mathbb X^*:=\{\xv:G(\xv)=0\}$ and $\xv^*\in\mathbb X^*$. The residual $G$ provides a local error bound at $\xv^*$ if, for some neighborhood $\mathbb N$ of $\xv^*$ and $c>0$,
\begin{equation}\label{Def-3.3a}
 c\,\dist(\xv,\mathbb X^*)\le\|G(\xv)\|,\qquad \xv\in\mathbb N.
\end{equation}
Here $\dist(\xv,\mathbb X^*):=\inf_{\yv\in\mathbb X^*}\|\xv-\yv\|$. 
\end{definition}

Definition \ref{LEBnew} applies the classical local error-bound condition for nonlinear equations, as used in \cite[Definition 1.1]{FY05}, to the component-gradient residual $G$. Equivalently, $G$ is metrically subregular at $(\xv^*,0)$ \cite{Kruger15}. It permits non-isolated roots. 

Suppose that the target $\boldsymbol\eta=(s_\nu t_0)_{\nu=1}^n$ is realizable, so that $\mathbb X^*\ne\varnothing$. Since $f_\nu(u)=q(s_\nu u)$ has its unique minimizer at $\eta_\nu=s_\nu t_0$, we have $F^*=nq(t_0)$. On a sufficiently small box around $\boldsymbol\eta$, continuity and positivity of the component curvatures give constants $0<m\le M_H$ such that $m\le f_\nu''(u)\le M_H$ on every segment joining $\eta_\nu$ to $\theta_\nu(\xv)$. We claim that
\begin{equation}\label{eq:eb-qg-connection}
 F(\xv)-F^*\asymp\|\thetav(\xv)-\boldsymbol\eta\|^2
                 \asymp\|G(\xv)\|^2,
\end{equation}
where $\asymp$ denotes two-sided bounds by positive constants independent of $\xv$ in this neighborhood.

To prove \eqref{eq:eb-qg-connection}, write $d_\nu:=\theta_\nu(\xv)-\eta_\nu$. Since $f_\nu'(\eta_\nu)=0$, Taylor's formula with integral remainder and the fundamental theorem of calculus give
\[
\begin{aligned}
 f_\nu(\eta_\nu+d_\nu)-f_\nu(\eta_\nu)
 &=d_\nu^2\int_0^1(1-t)f_\nu''(\eta_\nu+td_\nu)\,\mathrm dt,\\
 G_\nu(\xv)=f_\nu'(\eta_\nu+d_\nu)
 &=d_\nu\int_0^1 f_\nu''(\eta_\nu+td_\nu)\,\mathrm dt.
\end{aligned}
\]
The curvature bounds therefore imply, including when $d_\nu=0$,
\[
 \frac m2d_\nu^2\le f_\nu(\theta_\nu(\xv))-f_\nu(\eta_\nu)
 \le\frac{M_H}{2}d_\nu^2,
 \qquad m|d_\nu|\le|G_\nu(\xv)|\le M_H|d_\nu|.
\]
Summing the first inequalities and summing squares of the second yields
\[
\begin{aligned}
 \frac m2\|\thetav(\xv)-\boldsymbol\eta\|^2
 &\le F(\xv)-F^*\le\frac{M_H}{2}\|\thetav(\xv)-\boldsymbol\eta\|^2,\\
 m^2\|\thetav(\xv)-\boldsymbol\eta\|^2
 &\le\|G(\xv)\|^2\le M_H^2\|\thetav(\xv)-\boldsymbol\eta\|^2.
\end{aligned}
\]
These prove \eqref{eq:eb-qg-connection}. In particular,
\[
 \frac{m}{2M_H^2}\|G(\xv)\|^2
 \le F(\xv)-F^*\le\frac{M_H}{2m^2}\|G(\xv)\|^2.
\]
Thus the error bound in Definition \ref{LEBnew}, $\|G(\xv)\|\ge c\,\dist(\xv,\mathbb X^*)$, implies local quadratic growth with constant $mc^2/(2M_H^2)$. Conversely, if $F(\xv)-F^*\ge c_{\rm qg}\dist(\xv,\mathbb X^*)^2$ for some $c_{\rm qg}>0$, the same upper bound gives $c=\sqrt{2m^2c_{\rm qg}/M_H}$. This relates the component residual to local growth and nonisolated-minimum theory \cite{RB25}; both realizability and positive outer curvature are essential to this argument. The following theorem presents a standard sufficient criterion.

\begin{theorem}[A constant-rank sufficient condition]\label{thm:constant-rank-eb}
Let $\xv^*\in\mathbb X^*$. Suppose that, on a neighborhood $V$ of $\xv^*$,
\[
 H(\xv)\succeq m{E},\qquad \operatorname{rank}J(\xv)=r\ge1,
 \qquad{\sigma_r}(J(\xv))\ge\kappa>0,
\]
where $m,\kappa>0$ and $\sigma_r(J(\xv))$ is the smallest positive singular value of $J(\xv)$. Then, on a smaller neighborhood $\mathbb N$,
\begin{equation}\label{eq:constant-rank-eb}
 \|G(\xv)\|\ge\frac{m\kappa}{2}\dist(\xv,\mathbb X^*).
\end{equation}
\end{theorem}
\begin{proof}
We use the Euclidean open ball $B(\xv,a)=\{\yv:\|\yv-\xv\|<a\}$ and its closure $\overline B(\xv,a)$ in the factor space $\REAL^d$.
Since $\mathrm{D}G=HJ$, the derivative has constant rank $r$ and kernel $\ker J$. By the constant-rank theorem \cite[Theorem 9.32]{Rudin76}, there is an ambient open neighborhood $V_0\subset V$ such that $\mathcal Z:=G^{-1}(0)\cap V_0$ is a $C^1$ submanifold and $T_{\yv}\mathcal Z=\ker J(\yv)$. Choose $a>0$ with $\overline B(\xv^*,a)\subset V_0$ and
\[
 \|\mathrm{D}G(\mathbf u)-\mathrm{D}G(\mathbf v)\|\le m\kappa/2,
 \qquad\mathbf u,\mathbf v\in\overline B(\xv^*,a).
\]
For $\xv\in B(\xv^*,a/4)$, choose a nearest point $\bar\xv$ in the compact set $\mathcal Z\cap\overline B(\xv^*,a/2)$. Since $\xv^*$ belongs to this set,
\[
 \|\xv-\bar\xv\|\le\|\xv-\xv^*\|,\qquad
 \|\bar\xv-\xv^*\|\le2\|\xv-\xv^*\|<a/2.
\]
Thus $\bar\xv$ is an interior minimizer on the local manifold. With $\mathbf h=\xv-\bar\xv$, first-order optimality gives $\mathbf h\perp\ker J(\bar\xv)$, and hence
\[
 \|\mathrm{D}G(\bar\xv)\mathbf h\|\ge m\|J(\bar\xv)\mathbf h\|
 \ge m\kappa\|\mathbf h\|.
\]
Integration along the segment from $\bar\xv$ to $\xv$ \cite[Theorem 6.24]{Rudin76} gives
\[
 G(\xv)=\mathrm{D}G(\bar\xv)\mathbf h+
 \int_0^1[\mathrm{D}G(\bar\xv+t\mathbf h)-\mathrm{D}G(\bar\xv)]\mathbf h\,\mathrm d t.
\]
The integral has norm at most $m\kappa\|\mathbf h\|/2$. Therefore
\[
 \|G(\xv)\|\ge\tfrac12m\kappa\|\mathbf h\|
 \ge\tfrac12m\kappa\dist(\xv,\mathbb X^*),
\]
which proves the claim.
\end{proof}

The criterion permits rank-deficient Jacobians. If the rank at a root equals the maximal rank attained in an ambient neighborhood, a nonzero minor at that root remains nonzero locally, while the neighborhood rank bound excludes an increase. The Jacobian therefore has locally constant rank.

For rank-one tensors, monomial parameterization, exponent matrices, and sign constraints are established completion tools \cite[Section 2, especially Proposition 2.7]{KKKR17}. The following proposition applies Theorem \ref{thm:constant-rank-eb} to this structure. Its hypothesis is an explicit compatibility condition on the observed signs; realizability, local constant rank, and the residual error bound are conclusions.

\begin{proposition}[A rank-one error-bound certificate]\label{prop:rankone-pattern-eb}
Consider the model with $R=1$ and $\Omega\ne\varnothing$. Let $t_0>0$ be the unique minimizer of $q$, so $\theta_{ijk}=a_i b_j c_k$ and the target is $s_{ijk}t_0$. A factor coordinate is called active if it occurs in an observed entry. Let $A_\Omega$ have one row per observed index and one column per active coordinate; its row for $(i,j,k)$ has ones in the columns for $a_i,b_j,c_k$ and zeros elsewhere. Write $r=\operatorname{rank}_{\REAL}A_\Omega$.
	
Assume that the observed signs satisfy the compatibility condition
\[
(A_\Omega\mathbf z)_\nu\equiv(1-s_\nu)/2\pmod2,
\qquad\nu=1,\ldots,n,
\]	for some vector $\mathbf z$ with one binary component per active coordinate. Then:
\begin{enumerate}
\item The target is realizable, so $\mathbb X^*\ne\varnothing$. Every root has nonzero active coordinates.
\item At every $\xv^*\in\mathbb X^*$, the Jacobian has constant rank $r\ge1$ on a sufficiently small neighborhood. For any
\[
B>\max_{\ell\ {\rm active}}|x_\ell^*|,
\qquad 0<m\le\min_{t\in[t_0/2,\,3t_0/2]}q''(t),
\]
there is a neighborhood $\mathbb N$ of $\xv^*$ on which
\[
\|G(\xv)\|\ge
\frac{m t_0\sigma_r(A_\Omega)}{4B}
\operatorname{dist}(\xv,\mathbb X^*),\qquad\xv\in\mathbb N.	\]
\end{enumerate}
\end{proposition}
\begin{proof}
\textit{Realizability.} Let $\mathbf z$ solve the compatibility system and set each active coordinate to $(-1)^{z_\ell}t_0^{1/3}$. Every observed product is $s_\nu t_0$, hence $G=0$. Inactive coordinates are arbitrary. At any root, the nonzero observed products force all active coordinates to be nonzero. Conversely, the signs of these coordinates give a solution of the compatibility system. Thus sign compatibility is equivalent to realizability for these targets.
	
\textit{Verification of Theorem \ref{thm:constant-rank-eb}.} Fix a root $\xv^*$ and constants $B,m$ as stated. Such constants exist because the root is finite and $q''$ is continuous and positive on the compact interval. By continuity, choose a neighborhood on which all active coordinates remain nonzero, $|x_\ell|\le B$, and all signed margins belong to $[t_0/2,3t_0/2]$. Hence $H\succeq mE$.
	
Collect the active coordinates in $\xv_{\rm act}$ and the corresponding Jacobian columns in $J_{\rm act}$. Differentiating the products gives
\[
J_{\rm act}(\xv)=\operatorname{diag}(\thetav(\xv))\,
A_\Omega\operatorname{diag}(\xv_{\rm act})^{-1}.
\]
Both diagonal matrices are invertible, and the inactive-coordinate columns of $J$ vanish. Therefore $\operatorname{rank}J(\xv)=r$ throughout this neighborhood; $r\ge1$ since $\Omega\ne\varnothing$. Their smallest singular values are at least $t_0/2$ and $1/B$, respectively, so
\[
\sigma_r(J(\xv))\ge\frac{t_0}{2B}\sigma_r(A_\Omega).
\]
All hypotheses of Theorem \ref{thm:constant-rank-eb} have now been verified. Applying it with $\kappa=t_0\sigma_r(A_\Omega)/(2B)$ gives the asserted bound on a smaller neighborhood.
\end{proof}

Here $B$ and $m$ are choices of local bounding constants, rather than further restrictions on the observations. For example, the logistic part of $q''$ gives the admissible choice
\[
m=\frac{1}{\beta^2[2+2\cosh(3t_0/(2\beta))]}.
\]
The certificate separates loss curvature, target magnitude, observation-pattern spectrum, and factor scale. 

\begin{example}[{Partially observed rank-one latent model}]\label{ex:partial-rankone-eb}
Take $I=J=K=3$, $R=1$, and
\[
 \mathbf u=(1,-1,1)^T,\quad\mathbf v=(1,1,-1)^T,
 \quad\mathbf w=(1,-1,-1)^T.
\]
Set $s_{ijk}=u_iv_jw_k$ and $y_{ijk}=(s_{ijk}+1)/2$ on
\[
\begin{split}
 \Omega={}&\{(i,1,1):1\le i\le3\}\cup\{(1,j,1):1\le j\le3\}
             \cup\{(1,1,k):1\le k\le3\}\\
 &\quad{}\cup\{(2,2,1),(2,1,2),(1,2,2)\}.
\end{split}
\]
There are ten observations, nine factor variables, and seventeen missing entries. For ${\beta}=\tau=\lambda=1$, $\alpha=2$, the signed loss has $q'(1)=0$ and
\[
 q''(t)\ge m_0:=\frac{2e^2}{(1+e^2)^2}>0,\qquad0\le t\le2.
\]
{Indeed,
\[
 q''(t)=\frac1{2+2\cosh t}+\frac1{2+2\cosh(t-2)}.
\]
For $0\le t\le2$, both $|t|$ and $|t-2|$ are at most two, so each summand is at least $1/(2+2\cosh2)=e^2/(1+e^2)^2$.} The factors $(\mathbf u,\mathbf v,\mathbf w)$ realize every signed target. Within factorwise $\ell_\infty$ distance $1/4$ of them, signed entries lie in $(27/64,125/64)\subset(0,2)$, so $H\succeq m_0{E}$.

For this pattern, $A_\Omega$ has ten rows and nine columns. Its three observed fibers force a null vector to be constant within each mode, with the three constants summing to zero. Its nullspace therefore has dimension two, and $\operatorname{rank}A_\Omega=9-2=7$. Proposition \ref{prop:rankone-pattern-eb} now applies to the reference factors $(\mathbf u,\mathbf v,\mathbf w)$. After shrinking the neighborhood so that all signed margins lie in $[1/2,3/2]$, we may use $B=5/4$ and the curvature bound $m_0$ above. Hence
\[
 \|G(\xv)\|\ge\frac{m_0\sigma_7(A_\Omega)}5
                 \operatorname{dist}(\xv,\mathbb X^*).
\]
This is a concrete instance of the proposition with a row-rank-deficient Jacobian: its rank is seven although there are ten observations. 
\end{example}

\begin{example}[A triple-rank-two binary tensor]\label{ex:ranktwo-eb}
Take $I=J=K=R=2$, observe all eight entries, and prescribe the sign tensor $\mathcal S=(s_{ijk})$ by its two frontal slices
\[
 \mathcal S(:,:,1)=\begin{pmatrix}1&1\\1&-1\end{pmatrix},\qquad
 \mathcal S(:,:,2)=\begin{pmatrix}1&1\\1&1\end{pmatrix},\qquad
 y_{ijk}=(s_{ijk}+1)/2.
\]
Use the loss parameters of Example \ref{ex:partial-rankone-eb}, so the signed target is one. With $\delta_{pq}=1$ for $p=q$ and zero otherwise, set
\[
 b^*_{r_1jr_3}=\delta_{r_1,1}\delta_{r_3,j},\quad
 c^*_{r_1r_2k}=\delta_{r_1,1}\delta_{r_2,k},\quad
 a^*_{ir_2r_3}=s_{ir_3r_2}.
\]
The triple sum reduces to $\theta_{ijk}(\xv^*)=a^*_{ikj}=s_{ijk}$, hence $G(\xv^*)=0$. The first target slice has matrix rank two, which rules out $R=1$; the target therefore has triple rank exactly two.

With $\Bt,\Ct$ fixed, $\thetav=L(\Bt,\Ct)\vec\At$, and $L$ at the reference point is an $8\times8$ permutation matrix. By continuity, on a sufficiently small neighborhood $\sigma_{\min}(L)\ge1/2$ and all signed entries lie in $(0,2)$. Hence $JJ^T\succeq LL^T\succeq E/4$ and $H\succeq m_0E$. Theorem \ref{thm:constant-rank-eb} therefore applies with row rank eight. The local root manifold has dimension $24-8=16$. Here variations in one factor control all observed entries. Example \ref{ex:partial-rankone-eb} instead verifies the error bound in a partially observed, row-rank-deficient setting.
\end{example}

\section{Local quadratic convergence and model agreement}\label{Sec-4}

\begin{assumption}\label{ass:basic assumptions}
Fix $\xv^\circ\in\mathbb X^*$ and $U=B(\xv^\circ,\rho)$. Suppose that $\mathbb X^*$ is closed, the maps are $C^2$ on the relevant neighborhoods, and, on $U$,
\[
 m{E}\preceq H(\xv)\preceq M_H{E},\qquad\|J(\xv)\|\le M_J,
 \qquad c\,\dist(\xv,\mathbb X^*)\le\|G(\xv)\|,
\]
where $m,c>0$. Moreover, $H$ and $J$ are Lipschitz in factor space, with constants $L_H$ and $L_J$.
\end{assumption}
No constant-rank assumption is imposed here. Theorem \ref{thm:constant-rank-eb} is one verification route; direct corrections or other error-bound arguments are also admissible. {All roots in $U$ realize the same latent target $\boldsymbol\eta=\thetav(\xv^\circ)$, because each scalar loss has a unique minimizer. For each $\nu$, continuity on the connected ball $U$ makes $\theta_\nu(U)$ an interval containing $\eta_\nu$. The curvature bounds thus hold along the segment between $\eta_\nu$ and $\theta_\nu(\xv)$. The scalar integral argument used for \eqref{eq:eb-qg-connection} gives}
\begin{equation}\label{eq:latent-residual-equivalence}
 m\|\thetav(\xv)-\boldsymbol\eta\|\le\|G(\xv)\|
 \le M_H\|\thetav(\xv)-\boldsymbol\eta\|.
\end{equation}
Let $K_G=M_HM_J$ and $L_G=L_HM_J+M_HL_J$. Since $\mathrm{D}G=HJ$, for $\xv,\yv\in U$,
\begin{equation}\label{eq:residual-taylor-section4}
\begin{split}
 \|G(\yv)-G(\xv)\|&\le K_G\|\yv-\xv\|,\\
 \|G(\yv)-G(\xv)-H(\xv)J(\xv)(\yv-\xv)\|
 &\le\tfrac12L_G\|\yv-\xv\|^2.
\end{split}
\end{equation}

\subsection{A decomposition of the local error}\label{subsec:error-decomposition}
\begin{lemma}[Weighted step error decomposition]\label{lem:error-decomposition}
Under Assumption \ref{ass:basic assumptions}, take $\xv\in B(\xv^\circ,\rho/2)$ and a nearest root $\bar\xv$. Write $\mathbf h:=\xv-\bar\xv$, $e:=\|\mathbf h\|$, and $\boldsymbol\zeta:=\thetav(\xv)-\boldsymbol\eta$. Define two defects:
\[
 \underbrace{\mathbf a=\boldsymbol\zeta-H^{-1}G}_{\text{outer-loss defect}},\qquad
 \underbrace{\mathbf b=\boldsymbol\zeta-J\mathbf h}_{\text{factor-linearization defect}}.
\]
The Newton step for the outer problem in latent space is $-H^{-1}G$, so $H^{-1}G$ approximates the latent error. The defect $\mathbf a$ measures the error of this approximation; $\mathbf b$ measures the error of the factor linearization $J\mathbf h$. All unindexed matrices are evaluated at the current $\xv$. Define the operators
\[
 P=J(J^THJ)^\dagger J^TH,\qquad
 D_\mu=P-J(J^THJ+\mu {E})^{-1}J^TH,
\]
where $^\dagger$ denotes the Moore--Penrose pseudoinverse.
Here $P$ is the $H$-orthogonal projection onto $\operatorname{range}J$, and $D_\mu$ measures the change in this projection caused by damping. Both act in the observed latent-value space $\REAL^n$. Let $\mu>0$ and let $\sv_\mu$ solve $(J^THJ+\mu E)\sv_\mu=-J^TG$, with $\xv+\sv_\mu\in U$. Define
\[
 \mathbf R_\theta=\thetav(\xv+\sv_\mu)-\thetav(\xv)-J\sv_\mu
\]
as the nonlinear remainder of the exact step. Then
\begin{equation}\label{eq:three-source-decomposition}
 \thetav(\xv+\sv_\mu)-\boldsymbol\eta
 =P\mathbf a+({E}-P)\mathbf b+\mathbf R_\theta+D_\mu H^{-1}G.
\end{equation}
In the norm $\|\mathbf z\|_H=(\mathbf z^TH\mathbf z)^{1/2}$, a constant $C$, independent of $\xv$ and $\mu$ in the stated domain, satisfies
\begin{equation}\label{eq:component-bounds}
 \|\mathbf a\|_H+\|\mathbf b\|_H\le Ce^2,\qquad
 \|\mathbf R_\theta\|_H\le C\|\sv_\mu\|^2,
\end{equation}
\begin{equation}\label{eq:damping-basic-bound}
 \|D_\mu H^{-1}G\|_H\le\tfrac12\sqrt\mu\,e+Ce^2.
\end{equation}
\end{lemma}
\begin{proof}
The nearest root belongs to $U$: $e\le\|\xv-\xv^\circ\|$ implies $\|\bar\xv-\xv^\circ\|<\rho$. Taylor's theorem gives
$\|\mathbf b\|\le L_Je^2/2$ and $\|\mathbf R_\theta\|\le L_J\|\sv_\mu\|^2/2$.
Expanding $G(\bar\xv)=0$ at $\xv$ with displacement $-\mathbf h$, \eqref{eq:residual-taylor-section4} gives
$G=HJ\mathbf h+\mathbf r_G$, with $\|\mathbf r_G\|\le L_Ge^2/2$.
Consequently $\mathbf a=\mathbf b-H^{-1}\mathbf r_G$. Norm equivalence,
$\sqrt m\|\mathbf z\|\le\|\mathbf z\|_H\le\sqrt{M_H}\|\mathbf z\|$,
proves \eqref{eq:component-bounds}.

To see the projection property, set $W=H^{1/2}J$. A singular-value decomposition gives
$H^{1/2}PH^{-1/2}=WW^\dagger$, the Euclidean orthogonal projector onto $\operatorname{range}W$. Thus $PJ=J$, and both $P$ and $E-P$ are nonexpansive in the $H$-norm. The step equation yields
\[
 J\sv_\mu=-(P-D_\mu)H^{-1}G.
\]
Since $H^{-1}G=\boldsymbol\zeta-\mathbf a$ and $\boldsymbol\zeta=J\mathbf h+\mathbf b$, using $PJ=J$ yields
\[
\begin{aligned}
 \boldsymbol\zeta+J\sv_\mu
 &=\boldsymbol\zeta-P(\boldsymbol\zeta-\mathbf a)+D_\mu H^{-1}G\\
 &=({E}-P)(J\mathbf h+\mathbf b)+P\mathbf a+D_\mu H^{-1}G\\
 &=({E}-P)\mathbf b+P\mathbf a+D_\mu H^{-1}G.
\end{aligned}
\]
Adding $\mathbf R_\theta$ proves \eqref{eq:three-source-decomposition}.

On a left singular vector of $W$ with positive singular value $\sigma_i$, the transformed operator $H^{1/2}D_\mu H^{-1/2}$ acts by multiplication by $\mu/(\sigma_i^2+\mu)$; it is zero on the left nullspace. 
Its norm is at most one. Consequently, for every $\mathbf z$,
\[
 \|D_\mu\mathbf z\|_H
 =\|(H^{1/2}D_\mu H^{-1/2})H^{1/2}\mathbf z\|
 \le\|\mathbf z\|_H.
\]
Moreover,
\[
\|D_\mu J\mathbf h\|_H\le
\max_{\sigma\ge0}\frac{\mu\sigma}{\sigma^2+\mu}\,e
\le\tfrac12\sqrt\mu\,e.
\]
Finally, $H^{-1}G=J\mathbf h+\mathbf b-\mathbf a$, the triangle inequality, and \eqref{eq:component-bounds} give
\[
	\|D_\mu H^{-1}G\|_H
	\le\|D_\mu J\mathbf h\|_H+\|\mathbf b\|_H+\|\mathbf a\|_H
	\le\tfrac12\sqrt\mu\,e+Ce^2,
\]
which proves \eqref{eq:damping-basic-bound}.
\end{proof}

Equation \eqref{eq:three-source-decomposition} separates three error sources: $P\mathbf a$ comes from nonquadratic outer curvature; $({E}-P)\mathbf b+\mathbf R_\theta$ comes from nonlinear interactions among the factors; and $D_\mu H^{-1}G$ comes from damping. The projection isolates changes attainable to first order through the current factors. In particular, $({E}-P)\boldsymbol\zeta=({E}-P)\mathbf b$ is second order even for a rank-deficient Jacobian. A centered quadratic outer loss has $\mathbf a=0$ exactly.

\subsection{Quadratic convergence}\label{subsec:local-quadratic}
The next theorem provides the exact-step convergence guarantee used here. It follows the error-bound strategy of Fan--Yuan \cite{FY05}. We present the proof to make the damping range, rank treatment, and distinction between the component residual and the factor gradient explicit.
\begin{theorem}[Local quadratic convergence of the exact steps]\label{thm:Quad rate}\label{thm_local convergence}
Under Assumption \ref{ass:basic assumptions}, fix $1\le\delta\le2$. If $\xv_0$ is sufficiently close to $\xv^\circ$, the exact iteration \eqref{eq:LM method_system}--\eqref{eq:LM method} stays in $U$ and converges to some $\xv_\infty\in\mathbb X^*$. With $e_t=\dist(\xv_t,\mathbb X^*)$, there is $C>0$ such that, for all sufficiently large $t$,
\begin{equation}\label{eq:Quad rate}
 e_{t+1}\le Ce_t^2,\qquad
 \|\xv_{t+1}-\xv_\infty\|\le C\|\xv_t-\xv_\infty\|^2.
\end{equation}
Neither full row rank nor constant Jacobian rank is required.
\end{theorem}
\begin{proof}
The proof first establishes convergence using a coarse estimate, and then sharpens that estimate by separating the singular directions of the limiting weighted Jacobian. All constants below are independent of the iteration index. If an iterate belongs to $\mathbb X^*$, the method terminates. We therefore consider nonzero residuals.

\textit{Step 1: Nearest-point estimates.}
Let $\xv\in B(\xv^\circ,\rho/2)$ and choose a nearest root $\bar\xv\in\mathbb X^*$, which exists because $\mathbb X^*$ is nonempty and closed. As in Lemma \ref{lem:error-decomposition}, put $\mathbf h=\xv-\bar\xv$ and $e=\|\mathbf h\|$, and restrict attention to $e\le1$. Since $\xv^\circ$ is a root,
\[
 e\le\|\xv-\xv^\circ\|,\qquad
 \|\bar\xv-\xv^\circ\|\le e+\|\xv-\xv^\circ\|<\rho.
\]
Thus both points and their connecting segment lie in $U$. The error bound and the Lipschitz estimate in \eqref{eq:residual-taylor-section4}, with $G(\bar\xv)=0$, give
\[
 ce\le\|G\|=\|G(\xv)-G(\bar\xv)\|\le K_Ge.
\]
Taylor expansion at the current point, in the direction $-\mathbf h$, also gives
\[
 \|G-HJ\mathbf h\|\le\tfrac12L_Ge^2.
\]
Set $W=H^{1/2}J$, $\mathbf w=H^{-1/2}G$, and $A=L_G/(2\sqrt m)$. Multiplication by $H^{-1/2}$ has operator norm at most $1/\sqrt m$. Consequently, with $\mu=\|G\|^\delta$,
\begin{equation}\label{eq:weighted-taylor}
 c^\delta e^\delta\le\mu\le K_G^\delta e^\delta,
 \qquad\|\mathbf w-W\mathbf h\|\le Ae^2.
\end{equation}

\textit{Step 2: Step and residual estimates.}
Consider the strictly convex quadratic function
\[
 Q(\mathbf v)=\|W\mathbf v+\mathbf w\|^2+\mu\|\mathbf v\|^2.
\]
Its first-order condition is $(W^TW+\mu E)\mathbf v=-W^T\mathbf w$. Since $W^TW=J^THJ$ and $W^T\mathbf w=J^TG$, its unique minimizer is precisely the LM step $\sv$. Comparing this minimizer with $-\mathbf h$ and using \eqref{eq:weighted-taylor}, we obtain
\begin{equation}\label{eq:model-comparison}
 Q(\sv)\le Q(-\mathbf h)
 =\|\mathbf w-W\mathbf h\|^2+\mu e^2
 \le A^2e^4+\mu e^2.
\end{equation}
Both terms in $Q$ are nonnegative. Hence
\[
\begin{aligned}
 \|\sv\|^2
 &\le e^2+A^2e^4/\mu
 \le e^2+A^2c^{-\delta}e^{4-\delta}
 \le(1+A^2c^{-\delta})e^2,\\
 \|W\sv+\mathbf w\|^2
 &\le A^2e^4+K_G^\delta e^{2+\delta}
 \le(A^2+K_G^\delta)e^{2+\delta}.
\end{aligned}
\]
The final inequalities use $e\le1$, $\delta\le2$, and $4-\delta\ge2$. With $C_s=\sqrt{1+A^2c^{-\delta}}$, these bounds become
\begin{equation}\label{eq:bound of pk}
 \|\sv\|\le C_se,\qquad
 \|W\sv+\mathbf w\|\le\sqrt{A^2+K_G^\delta}\,e^{1+\delta/2}.
\end{equation}
Whenever $\xv+\sv\in U$, Taylor expansion along the actual step gives
\[
 \|G(\xv+\sv)\|
 \le\|G+HJ\sv\|+\tfrac12L_G\|\sv\|^2
 \le\sqrt{M_H}\|\mathbf w+W\sv\|
       +\tfrac12L_GC_s^2e^2.
\]
Here $G+HJ\sv=H^{1/2}(\mathbf w+W\sv)$. Applying the error bound at $\xv+\sv$, and using $e^2\le e^{1+\delta/2}$, yields
\begin{equation}\label{eq:bound of dist zk}
 e_+\le C_pe^{1+\delta/2},\qquad
 C_p=\frac{\sqrt{M_H}\sqrt{A^2+K_G^\delta}+L_GC_s^2/2}{c},
\end{equation}
where $e_+=\dist(\xv+\sv,\mathbb X^*)$. We next verify that all trial points stay in $U$, so this estimate applies throughout the iteration.

\textit{Step 3: Local invariance and convergence.}
Choose
\[
 0<a<\min\{\rho/[4(1+C_s)],(2C_p)^{-2/\delta},1\},
 \qquad\|\xv_0-\xv^\circ\|<a.
\]
We prove inductively that $\xv_t\in B(\xv^\circ,\rho/2)$ and $e_t\le2^{-t}e_0$. Suppose these properties hold up to index $t$. The step bound is available before requiring the next point to lie in $U$, and gives
\[
 \|\xv_{t+1}-\xv^\circ\|
 \le\|\xv_0-\xv^\circ\|+\sum_{i=0}^t\|\sv_i\|
 \le a+C_s\sum_{i=0}^t e_i
 <a+2C_sa<\rho/2.
\]
Thus the next point lies in $U$. Since $e_t\le e_0<a$ and $C_pa^{\delta/2}<1/2$, \eqref{eq:bound of dist zk} now implies
\[
 e_{t+1}\le C_pe_t^{\delta/2}e_t\le\tfrac12e_t.
\]
This closes the induction. In particular, $\sum_t\|\sv_t\|\le C_s\sum_t e_t<\infty$, so $\xv_t$ is Cauchy and converges to a point $\xv_\infty$ inside $U$. Since $e_t\to0$ and $\mathbb X^*$ is closed, $\xv_\infty\in\mathbb X^*$.

The same geometric decrease holds from every index $t$ onward: $e_{t+j}\le2^{-j}e_t$. Therefore
\[
 \|\xv_t-\xv_\infty\|
 \le\sum_{i=t}^\infty\|\sv_i\|
 \le C_s\sum_{j=0}^\infty e_{t+j}\le2C_se_t.
\]
Because $\xv_\infty$ is a root, the reverse inequality $e_t\le\|\xv_t-\xv_\infty\|$ is immediate. Writing $\varepsilon_t$ for the latter distance, we have
\begin{equation}\label{eq:tail-distance-equivalence}
 e_t\le\varepsilon_t:=\|\xv_t-\xv_\infty\|\le2C_se_t.
\end{equation}
For $\delta=2$, \eqref{eq:bound of dist zk} already gives $e_{t+1}\le C_pe_t^2$. Moreover,
\[
 \varepsilon_{t+1}\le2C_se_{t+1}
 \le2C_sC_pe_t^2\le2C_sC_p\varepsilon_t^2.
\]
Thus both quadratic conclusions hold in this case. The remaining steps sharpen the coarse estimate for $1\le\delta<2$.

\textit{Step 4: Spectral decomposition at the limit.}
Let $W_\infty=H(\xv_\infty)^{1/2}J(\xv_\infty)$ and $r=\operatorname{rank}W_\infty$. Since $H$ is diagonal with entries at least $m$, the scalar identity
\[
 |\sqrt u-\sqrt v|=\frac{|u-v|}{\sqrt u+\sqrt v}
 \le\frac{|u-v|}{2\sqrt m}\qquad(u,v\ge m)
\]
implies $\|H_t^{1/2}-H_\infty^{1/2}\|\le L_H\varepsilon_t/(2\sqrt m)$. Expanding the difference of the matrix products gives
\[
\begin{aligned}
 \|W_t-W_\infty\|
 &\le\|H_t^{1/2}(J_t-J_\infty)\|
       +\|(H_t^{1/2}-H_\infty^{1/2})J_\infty\|\\
 &\le L_W\varepsilon_t,\qquad
 L_W=\sqrt{M_H}L_J+M_JL_H/(2\sqrt m).
\end{aligned}
\]
First suppose $r>0$. Set $b=\sigma_r(W_\infty)/2>0$, and let $\Pi_t$ be the orthogonal projector onto the first $r$ left singular vectors of $W_t$. Singular-value perturbation bounds imply, for all sufficiently large $t$,
\[
 \sigma_r(W_t)\ge\sigma_r(W_\infty)-\|W_t-W_\infty\|\ge b,
 \qquad\|(E-\Pi_t)W_t\|\le L_W\varepsilon_t.
\]
For the second inequality, if $r<\min\{n,d\}$, use $\sigma_{r+1}(W_\infty)=0$ and
$\|(E-\Pi_t)W_t\|=\sigma_{r+1}(W_t)\le\|W_t-W_\infty\|$; if $r=\min\{n,d\}$, the product $(E-\Pi_t)W_t$ is zero. Thus singular values outside these first $r$ directions may vanish, but their size is at most of order $\varepsilon_t$.

The complementary component of $\mathbf w_t$ satisfies a second-order estimate. Choose the nearest-root displacement $\mathbf h_t=\xv_t-\bar\xv_t$, so $\|\mathbf h_t\|=e_t$. By \eqref{eq:weighted-taylor},
\[
\begin{aligned}
 \|(E-\Pi_t)\mathbf w_t\|
 &\le\|(E-\Pi_t)(\mathbf w_t-W_t\mathbf h_t)\|
       +\|(E-\Pi_t)W_t\mathbf h_t\|\\
 &\le Ae_t^2+L_W\varepsilon_te_t
 \le(A+L_W)\varepsilon_t^2.
\end{aligned}
\]
This estimate uses the residual's approximation by $W_t\mathbf h_t$, rather than a uniform lower bound for the singular values outside the limiting positive spectrum.

\textit{Step 5: Refined regularization-error estimate.}
As shown in Lemma \ref{lem:error-decomposition}, in weighted coordinates the operator $D_{\mu_t}$ has multiplier $\mu_t/(\sigma_i(W_t)^2+\mu_t)$ on each positive left singular direction and multiplier zero on the left nullspace. On the first $r$ directions, this multiplier is at most $\mu_t/b^2$; on the remaining directions, it is at most one. Splitting $\mathbf w_t$ with $\Pi_t$ therefore gives
\[
 \|D_{\mu_t}H_t^{-1}G_t\|_{H_t}
 \le\frac{\mu_t}{b^2}\|\Pi_t\mathbf w_t\|
       +\|(E-\Pi_t)\mathbf w_t\|
 \le\frac{K_G^{1+\delta}}{b^2\sqrt m}\,
       \varepsilon_t^{1+\delta}
       +(A+L_W)\varepsilon_t^2.
\]
Here we used $\mu_t\le K_G^\delta e_t^\delta$ and $\|\mathbf w_t\|\le\|G_t\|/\sqrt m\le K_Ge_t/\sqrt m$, together with $e_t\le\varepsilon_t$. Since $\delta\ge1$ and eventually $\varepsilon_t\le1$, the whole right-hand side is bounded by a constant times $\varepsilon_t^2$.

If $r=0$, then $W_\infty=0$. Take $\Pi_t=0$ and omit $b$. Directly,
\[
 \|\mathbf w_t\|\le\|\mathbf w_t-W_t\mathbf h_t\|
              +\|W_t\|e_t
 \le Ae_t^2+L_W\varepsilon_te_t
 \le(A+L_W)\varepsilon_t^2.
\]
All damping multipliers are at most one, so the same quadratic damping bound follows. No constant-rank hypothesis has been used in either case.

\textit{Step 6: Quadratic convergence rates.}
In \eqref{eq:three-source-decomposition}, the first two defects have $H_t$-norm bounded by $Ce_t^2$ in total, because $P_t$ and $E-P_t$ are contractions in that norm. The nonlinear step remainder satisfies
\[
 \|\mathbf R_{\theta,t}\|_{H_t}
 \le C\|\sv_t\|^2\le CC_s^2e_t^2.
\]
Combining these bounds with Step 5 shows that, for some fixed $K>0$ and all sufficiently large $t$,
\[
 \|\thetav(\xv_{t+1})-\boldsymbol\eta\|_{H_t}
 \le K\varepsilon_t^2.
\]
The error bound, \eqref{eq:latent-residual-equivalence}, and norm equivalence now give
\[
 e_{t+1}
 \le\frac1c\|G(\xv_{t+1})\|
 \le\frac{M_H}{c}\|\thetav(\xv_{t+1})-\boldsymbol\eta\|
 \le\frac{M_HK}{c\sqrt m}\varepsilon_t^2.
\]
Finally, applying \eqref{eq:tail-distance-equivalence} on the appropriate side yields
\[
 e_{t+1}\le\frac{4M_HKC_s^2}{c\sqrt m}\,e_t^2,
 \qquad
 \varepsilon_{t+1}\le2C_se_{t+1}
 \le\frac{2M_HKC_s}{c\sqrt m}\,\varepsilon_t^2.
\]
Increasing a single constant $C$ to cover these two bounds proves \eqref{eq:Quad rate} for $1\le\delta<2$, completing the proof.
\end{proof}

The comparison \eqref{eq:weighted-taylor} and the singular-value estimates connect the established error-bound analysis \cite{FY05,WHZW24} to the present system. The tail estimate \eqref{eq:tail-distance-equivalence} then converts the bounds relative to the limit point into the asserted bounds relative to the root set. The ratio $M_H/m$ controls outer-weight distortion, not the whole linear-system condition number. For $W=H^{1/2}J$, $m{\sigma_i}(J)^2\le{\sigma_i}(W)^2\le M_H{\sigma_i}(J)^2$. If $\ker J\ne\{0\}$, then
\[
 \kappa_2(J^THJ+\mu {E})=1+\lambda_{\max}(J^THJ)/\mu.
\]
Thus even $M_H/m=1$ does not prevent ill-conditioning as damping vanishes.

\begin{corollary}[Acceptance of exact full LM steps]\label{cor:TR-full}
Under Assumption~\ref{ass:basic assumptions}, use $H_t^{\rm tr}=H_t$ and exact solves in Algorithm~\ref{alg:TRCGLM}. For each fixed $\Delta_0>0$, if $\xv_0$ is sufficiently close to $\xv^\circ$, the algorithm accepts the unrestricted LM step at every nonterminal iteration. With zero tolerance and no iteration limit, both quadratic rates of Theorem~\ref{thm:Quad rate} hold.
\end{corollary}
\begin{proof}
Use $W$, $\mathbf w$, $Q$, and $e$ from the proof of Theorem~\ref{thm:Quad rate}. For the exact unrestricted step $\sv$, \eqref{eq:model-comparison} gives $Q(\sv)=\mathcal O(e^{2+\delta})$. Since $\|\mathbf w\|^2\ge(c^2/M_H)e^2$, the predicted reduction
\[
 -\phi_t(\sv)=\tfrac12\bigl(\|\mathbf w\|^2-Q(\sv)\bigr)
\]
is positive and bounded below by a positive multiple of $e^2$ at sufficiently close nonroot points. Also $F(\xv)-F^*=\tfrac12\|\mathbf w\|^2+\mathcal O(e^3)$ by Taylor expansion and the outer-loss defect bound. The next-residual estimate in the same proof and \eqref{eq:eb-qg-connection} give $F(\xv+\sv)-F^*=\mathcal O(e^{2+\delta})$. Hence the reduction ratio at this trial is $1+\mathcal O(e)$.

The unrestricted iteration stays close to $\xv^\circ$ by Theorem~\ref{thm:Quad rate} and satisfies $\|\sv\|\le C_se$. Choose the initial neighborhood small enough that its steps have norm less than $\Delta_0$ and ratio above $3/4$. Each such step is the exact constrained minimizer, is accepted, and leaves the radius unchanged. Induction shows that the trust-region iterates coincide with the unrestricted iterates.
\end{proof}

\subsection{Improved model agreement for a symmetric logistic--LSE loss}\label{subsec:fullstep-agreement}
For common loss parameters, the root target is $\eta_\nu=s_\nu t_0$. {Write $\boldsymbol\zeta=\thetav(\xv)-\boldsymbol\eta$ for the latent error.} If $q'''(t_0)=0$, scalar Taylor expansion gives
\begin{equation}\label{eq:symmetric-loss-expansion}
 G=h_0\boldsymbol\zeta+\mathcal O(\|\boldsymbol\zeta\|^3),\qquad
 H=h_0{E}+\mathcal O(\|\boldsymbol\zeta\|^2),\qquad h_0=q''(t_0)>0.
\end{equation}
Thus the outer defect $\mathbf a=\boldsymbol\zeta-H^{-1}G$ is cubic, whereas mixed-factor terms are generally quadratic. For example, ${\beta}=\tau=\lambda=1$ gives $q(\alpha-t)=q(t)$, $t_0=\alpha/2$, and $q'''(t_0)=0$ for every $\alpha>0$.

Actual-to-predicted decrease ratios give a standard measure of model accuracy \cite[Section 4.4]{NW06}. For an exact step with $\delta=2$, set
\begin{equation}\label{eq:ratio-comparison-definitions}
 {\phi_{\xv}}(\sv)=G^TJ\sv+\tfrac12\sv^T(J^THJ+\mu {E})\sv,
 \quad P_{\xv}=-{\phi_{\xv}}(\sv),\quad
 \mathcal E_{\xv}=F(\xv+\sv)-F(\xv)-{\phi_{\xv}}(\sv).
\end{equation}
Here $\phi_{\xv}$ is the same damped quadratic step model as $\phi_t$, $P_{\xv}$ is its predicted decrease, and $\mathcal E_{\xv}$ is the objective-model discrepancy. Thus $\rho_{\xv}=1-\mathcal E_{\xv}/P_{\xv}$ is the actual-to-predicted decrease ratio. For analytic losses without cancellation, Taylor expansion gives $\mathcal E_{\xv}=\mathcal O(e^3)$. The predicted decrease is of order $e^2$, as shown below, so $|1-\rho_{\xv}|=\mathcal O(e)$. The next result improves both orders for the stated nonquadratic loss. The improvement uses two cancellations: the third outer derivative vanishes at the target, and the exact local step makes $G+HJs$ second order. This result concerns the order of the objective-model discrepancy and the decrease ratio; the iterate convergence order remains quadratic.

\begin{proposition}[Fourth-order model error and trust-region acceptance]\label{prop:fullstep-agreement}
Under Assumption \ref{ass:basic assumptions}, suppose the common logistic--LSE loss satisfies $q'''(t_0)=0$. For $\delta=2$, use the exact original-curvature step. On a sufficiently small neighborhood, for $e=\dist(\xv,\mathbb X^*)>0$,
\begin{equation}\label{eq:predicted-decrease-two-sided}
 \frac{c^2}{4M_H}e^2\le P_{\xv}\le\frac{K_G^2}{2m}e^2,
\end{equation}
and, for finite uniform constants,
\begin{equation}\label{eq:fourth-order-model-error}
 |\mathcal E_{\xv}|\le C_{\rm mod}e^4,\qquad
 |1-\rho_{\xv}|\le C_{\rm ratio}e^2,
 \quad \rho_{\xv}=1-\mathcal E_{\xv}/P_{\xv}.
\end{equation}
In the trust-region method with the original curvature, the exact unrestricted step is accepted if it is feasible and $C_{\rm ratio}e^2\le1-\eta_{\rm a}$.
\end{proposition}
\begin{proof}
Use $W=H^{1/2}J$, $\mathbf w=H^{-1/2}G$, and the quadratic function $Q$ from the proof of Theorem \ref{thm:Quad rate}. For $\delta=2$, \eqref{eq:model-comparison} gives $Q(\sv)\le(A^2+K_G^2)e^4$. Expanding $Q$ and using the definition of $\phi_{\xv}$ gives
\[
 Q(\sv)=\|\mathbf w\|^2+2\phi_{\xv}(\sv),\qquad
 P_{\xv}=\tfrac12(\|\mathbf w\|^2-Q(\sv)).
\]
Moreover, $ce\le\|G\|\le K_Ge$ and $mE\preceq H\preceq M_HE$ imply
\[
 \frac{c^2}{M_H}e^2\le\|\mathbf w\|^2
 \le\frac{K_G^2}{m}e^2.
\]
For sufficiently small $e$, $(A^2+K_G^2)e^4\le c^2e^2/(2M_H)$. Substituting these bounds into the identity for $P_{\xv}$ proves both inequalities in \eqref{eq:predicted-decrease-two-sided}.

Write $\mathbf u=J\sv$ and retain the nonlinear remainder $\mathbf R_\theta$ from Lemma \ref{lem:error-decomposition}, so that
\[
 \thetav(\xv+\sv)-\thetav(\xv)=\mathbf u+\mathbf R_\theta.
\]
The step bound gives $\|\sv\|=\mathcal O(e)$, hence $\|\mathbf u\|=\mathcal O(e)$ and $\|\mathbf R_\theta\|=\mathcal O(e^2)$. For $\delta=2$, \eqref{eq:bound of pk} also yields
\begin{equation}\label{eq:linearized-gradient-cancellation}
 \|G+H\mathbf u\|=\|H^{1/2}(\mathbf w+W\sv)\|=\mathcal O(e^2).
\end{equation}
Since $f_\nu'''(\eta_\nu)=s_\nu q'''(t_0)=0$ and the fourth derivatives are locally bounded, each third derivative along the latent Taylor segment is $\mathcal O(e)$. Taylor expansion of the component losses and summation therefore give
\[
 F(\xv+\sv)-F(\xv)
 =G^T(\mathbf u+\mathbf R_\theta)
 +\tfrac12(\mathbf u+\mathbf R_\theta)^TH
                  (\mathbf u+\mathbf R_\theta)+R_f,
 \qquad |R_f|=\mathcal O(e^4).
\]
Indeed, the latent distance to the target and the latent step are both $\mathcal O(e)$, so the third-order remainder is an $\mathcal O(e)$ coefficient times an $\mathcal O(e^3)$ term. Subtracting the damped quadratic model and collecting terms gives
\begin{equation}\label{eq:model-error-cancellation}
 \mathcal E_{\xv}=(G+H\mathbf u)^T\mathbf R_\theta
       +\tfrac12\mathbf R_\theta^TH\mathbf R_\theta+R_f
       -\tfrac12\mu\|\sv\|^2.
\end{equation}
The first term is $\mathcal O(e^4)$ by \eqref{eq:linearized-gradient-cancellation}; the second and third terms are $\mathcal O(e^4)$ by the remainder estimates. The final term is $\mathcal O(e^4)$ because $\mu=\mathcal O(e^2)$ and $\|\sv\|=\mathcal O(e)$. Thus $|\mathcal E_{\xv}|\le C_{\rm mod}e^4$. Dividing by $P_{\xv}\ge c^2e^2/(4M_H)$ gives $|1-\rho_{\xv}|\le(4M_HC_{\rm mod}/c^2)e^2$.

For a feasible unrestricted step, strict convexity makes it the exact minimizer of \eqref{eq:LM method_trust region}. Its ratio satisfies $\rho_{\xv}\ge1-C_{\rm ratio}e^2\ge\eta_{\rm a}$, so Algorithm~\ref{alg:TRCGLM} accepts it.
\end{proof}
Without $q'''(t_0)=0$, the same calculation gives
\[
 \mathcal E_{\xv}=\tfrac16\sum_\nu f_\nu'''(\eta_\nu)(J\sv)_\nu^3
                 +\mathcal O(e^4),
\]
which identifies the possible third-order leading term. As an illustrative scalar instance, take $\theta(x)=x$, one positive label, ${\beta}=\tau=1$, $\lambda=2$, and $\alpha=\log4$. Then $t_0=\log((\sqrt{33}-1)/4)$ and $q'''(t_0)>0$. For $x=t_0+z$, Taylor expansion gives
\[
 \rho_x-1=\frac{q'''(t_0)}{3q''(t_0)}z+\mathcal O(z^2).
\]
Thus the first-order ratio term need not vanish. Quadratic losses already have this cancellation; the proposition identifies a nonquadratic logistic--LSE family with the same property. The argument uses local fourth-order outer smoothness and second-order inner remainders and is not exclusive to triple parameterization. The proposition therefore establishes a model-accuracy refinement while preserving the quadratic order of the iterates.

\section{Numerical experiments}\label{Sec-7}

A synthetic example illustrates the local convergence rate of the exact LM iteration. We then assess prediction and computational performance of the trust-region implementation on HCP and Enron.

\subsection{Local quadratic convergence}\label{subsec:review-local}
We use a synthetic binary tensor to illustrate the local convergence of the exact LM iteration. The tensor size is $I\times J\times K=10\times8\times5$ and the triple rank is $R=2$, giving 400 observations and 92 factor variables. To specify the reference factors, let $p_\ell=(\ell-1)\bmod2$ and $q_\ell=\lfloor(\ell-1)/2\rfloor\bmod2$. For $r_1,r_2,r_3\in\{1,2\}$, set
\[
\begin{aligned}
	a^*_{ir_2r_3}&=\mathbf1_{\{r_2=p_i+1\}}(-1)^{q_i(r_3-1)},\\
	b^*_{r_1jr_3}&=\mathbf1_{\{r_3=p_j+1\}}(-1)^{q_j(r_1-1)},\\
	c^*_{r_1r_2k}&=\mathbf1_{\{r_1=p_k+1\}}(-1)^{q_k(r_2-1)},
\end{aligned}
\]
where $\mathbf1_{\{\cdot\}}$ equals one when its condition holds and zero otherwise. The resulting latent entries are
\[
\theta_{ijk}(\xv^*)=(-1)^{q_i p_j+q_j p_k+q_k p_i}\in\{-1,1\}.
\]
We observe all entries and set $y_{ijk}=[\theta_{ijk}(\xv^*)+1]/2$. This gives 248 labels equal to one and 152 equal to zero. The loss parameters are $\beta=\tau=\lambda=1$ and $\alpha=2$. Since $q'(1)=0$ and $q''>0$, the unique scalar minimizer is $t_0=1$. Consequently, the target vector is $\boldsymbol\eta=\thetav(\xv^*)$ and the reference factors satisfy $G(\xv^*)=0$.

We perturb the reference factors by $\xv_0=\xv^*+0.2\mathbf u$, where the entries of $\mathbf u\in\REAL^{92}$ are independent and uniformly distributed on $[0,1)$. From this common initial point, we run five iterations for each damping exponent $\delta\in\{1,1.5,2\}$. At each iteration, the MATLAB code evaluates $J_t$, $G_t$, and the original weight $H_t$, solves \eqref{eq:LM method_system} directly with $\mu_t=\|G_t\|^\delta$, and takes the full step $\xv_{t+1}=\xv_t+\sv_t$.

\begin{figure}[t]\centering
	\includegraphics[width=0.82\linewidth]{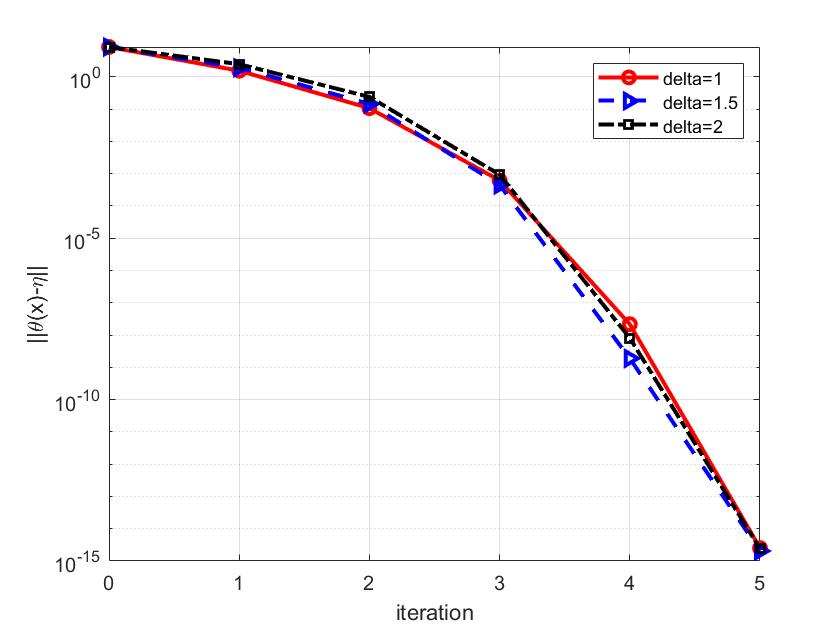}
	\caption{Residual $r_t=\|\thetav(\xv_t)-\boldsymbol\eta\|$ versus iteration for the $10\times8\times5$ binary tensor with $R=2$. The three runs use the same initial point and $\delta=1,1.5,2$.}\label{fig:local-rates}
\end{figure}

We record $r_t=\|\thetav(\xv_t)-\boldsymbol\eta\|$, the error in the fitted latent tensor. Under the local error bound and curvature assumptions, the residual is locally equivalent to the distance from the solution set; quadratic convergence therefore implies $r_{t+1}=\mathcal O(r_t^2)$. Figure~\ref{fig:local-rates} shows the MATLAB residual histories on a logarithmic scale. For each of the three exponents, the residual decreases to approximately $10^{-15}$ within five iterations. This behavior supports the local quadratic convergence in Theorem~\ref{thm:Quad rate}. The three damping exponents produce comparable local residual histories.

\subsection{Binary completion on real data}\label{subsec:review-data}
The HCP tensor has dimensions $68\times68\times212$ and records the presence or absence of structural connections between brain regions for each subject \cite{WZD19,WL20}. Its first two modes are symmetric. The Enron tensor is a $191\times1536\times40$ subset of the CMU email corpus,\footnote{\path"https://www.cs.cmu.edu/~enron/".} with a positive entry for a sender--recipient pair for which the sender sent at least one message to the recipient in the indicated period. About $0.35\%$ of its entries are positive. In the reported experiments, $80\%$ of the entries are used for training and the remaining $20\%$ are held out for evaluation.

The comparison includes four solvers: trust-region LM, majorization--minimization (MM) \cite{WL20,ZL21}, nonlinear conjugate gradients (NCG), and limited-memory BFGS (L-BFGS). For each solver, the tables include the canonical polyadic (CP) and triple (TTD) representations. The triple sizes selected by a BIC-based search over $1,\ldots,20$ are $R=6$ for HCP and $R=3$ for Enron. The matched CP ranks, $36$ and $9$, respectively, give the same number of factor entries as the corresponding triple models. For the penalized model, $\beta=0.1$, $\alpha=1$, $\tau=30$, and $\lambda=1$. The experimental LM model uses the diagonal curvature floor $h_{\rm floor}=0.1$ and truncated CG as described in Section~\ref{subsec:trust-implementation}. The reported LM runs use a relative gradient tolerance of $10^{-8}$ and a limit of 30 outer iterations.

Let $\Omega_{\mathrm{test}}$ be the set of held-out indices, and let $\widehat y_{ijk}$ be the binary prediction obtained by thresholding the fitted logistic probability at $0.5$. We define the root mean squared error (RMSE) and the misclassification error rate (MER) by
\[
\begin{aligned}
	\mathrm{RMSE}
	&:= \left(\frac{1}{|\Omega_{\mathrm{test}}|}
	\sum_{(i,j,k)\in\Omega_{\mathrm{test}}}(\widehat y_{ijk}-y_{ijk})^2\right)^{1/2},\\
	\mathrm{MER}
	&:= \frac{1}{|\Omega_{\mathrm{test}}|}
	\#\Big\{(i,j,k)\in\Omega_{\mathrm{test}}:\widehat y_{ijk}\ne y_{ijk}\Big\}.
\end{aligned}
\]
We also report the area under the receiver operating characteristic curve (AUC), computed from the fitted probabilities before thresholding.

\begin{table}[!htbp]
\centering\small
\caption{Held-out prediction and computational cost on HCP. MER is reported as a fraction.}\label{tb:HCP}
\begin{tabular}{llrrrrr}\hline
Representation & Method & RMSE & MER & AUC & Time (s) & Iterations \\ \hline
CP $(R=36)$ & MM & 0.28736 & 0.08257 & 97.30\% & 2.2003 & 8 \\
 & NCG & 0.28445 & 0.08091 & 97.48\% & 46.1043 & 100 \\
 & L-BFGS & 0.28323 & 0.08022 & 97.51\% & 32.8580 & 100 \\
 & LM & 0.27715 & 0.07680 & 97.68\% & 40.7107 & 21 \\
\hline
TTD $(R=6)$ & MM & 0.28521 & 0.08134 & 97.32\% & 2.1473 & 12 \\
 & NCG & 0.28041 & 0.07863 & 97.63\% & 40.384 & 100 \\
 & L-BFGS & 0.27979 & 0.07828 & 97.66\% & 28.769 & 100 \\
 & LM & 0.27625 & 0.07631 & 97.80\% & 43.871 & 30 \\
\hline
\end{tabular}
\end{table}

\begin{table}[!htbp]
\centering\small
\caption{Held-out prediction and computational cost on Enron. MER is reported as a fraction.}\label{tb:Enron}
\begin{tabular}{llrrrrr}\hline
Representation & Method & RMSE & MER & AUC & Time (s) & Iterations \\ \hline
CP $(R=9)$ & MM & $6.4869\times10^{-2}$ & $4.2079\times 10^{-3}$ & $97.22\%$ & $22.275$ & 17 \\
 & NCG & $6.2465\times10^{-2}$ & $3.9019\times 10^{-3}$ & $98.16\%$ & $785.47$ & 100 \\
 & L-BFGS & $6.2220\times10^{-2}$ & $3.8713\times 10^{-3}$ & $97.57\%$ & $477.43$ & 100 \\
 & LM & $3.8194\times10^{-2}$ & $1.4588\times 10^{-3}$ & $99.81\%$ & $647.23$ & 24 \\
\hline
TTD $(R=3)$ & MM & $6.0978\times10^{-2}$ & $3.7183\times 10^{-3}$ & $99.24\%$ & $52.061$ & 30 \\
 & NCG & $5.6776\times10^{-2}$ & $3.2235\times 10^{-3}$ & $99.48\%$ & $374.45$ & 100 \\
 & L-BFGS & $5.5366\times10^{-2}$ & $3.0654\times 10^{-3}$ & $99.36\%$ & $588.17$ & 100 \\
 & LM & $3.4399\times10^{-2}$ & $1.1833\times 10^{-3}$ & $99.88\%$ & $632.16$ & 26 \\
\hline
\end{tabular}
\end{table}

On HCP, TTD--LM achieves the lowest MER ($7.631\%$) and the highest AUC ($97.80\%$) in Table~\ref{tb:HCP}. It takes 30 outer iterations and $43.87$ seconds, compared with 100 iterations and $40.38$ and $28.77$ seconds for TTD--NCG and TTD--L-BFGS, respectively. On Enron (Table~\ref{tb:Enron}), TTD--LM reduces the MER from $0.3065\%$ for TTD--L-BFGS to $0.1183\%$ and increases the AUC from $99.36\%$ to $99.88\%$. It takes 26 outer iterations and $632.16$ seconds, compared with 100 iterations and $588.17$ seconds for TTD--L-BFGS. In the Enron LM runs, TTD also gives lower prediction error than CP at slightly lower computational cost ($632.16$ versus $647.23$ seconds).

\begin{figure}[H]
\centering
\begin{tabular}{@{}c@{\hspace{1em}}c@{}}
\includegraphics[width=.39\linewidth]{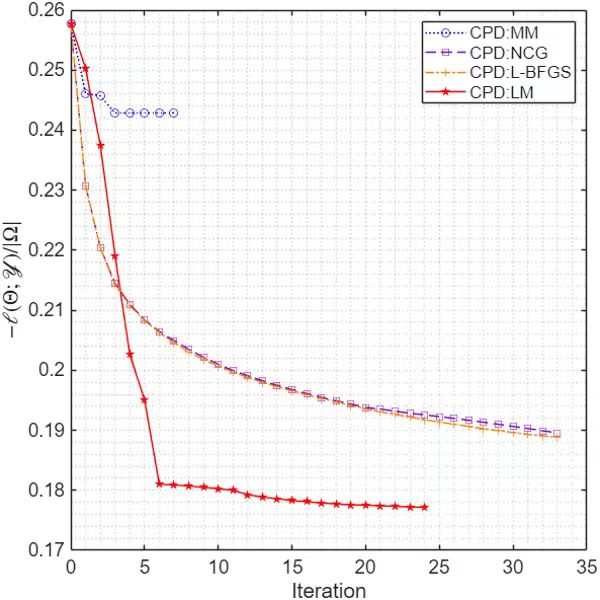} &
\includegraphics[width=.39\linewidth]{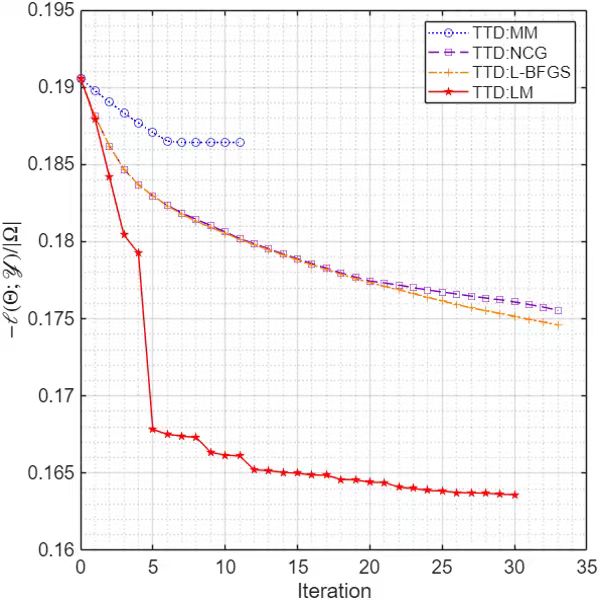} \\
\small HCP, CP & \small HCP, TTD \\[.5em]
\includegraphics[width=.39\linewidth]{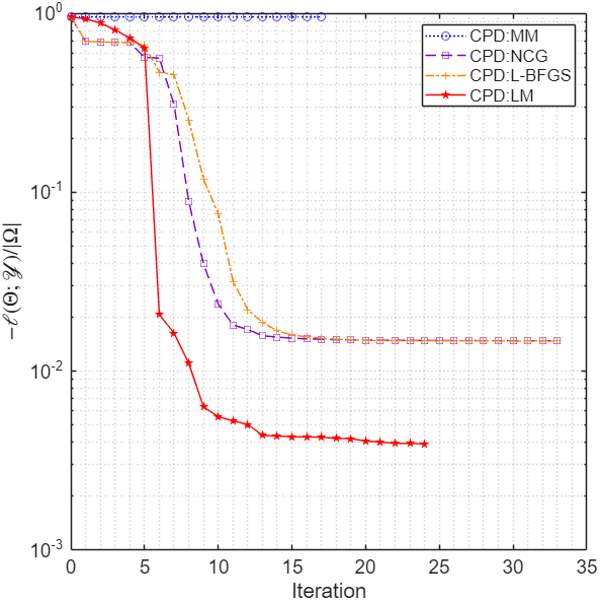} &
\includegraphics[width=.39\linewidth]{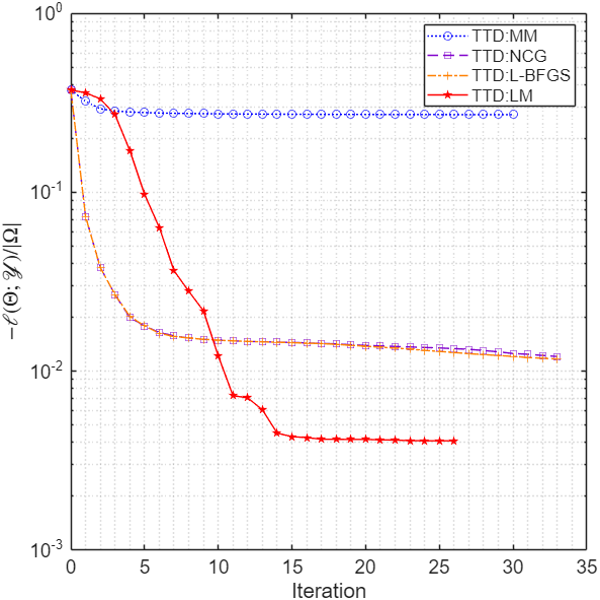} \\
\small Enron, CP & \small Enron, TTD
\end{tabular}
\caption{Average training negative log-likelihood versus outer iteration for HCP (top) and Enron (bottom). The Enron panels use a logarithmic vertical scale.}
\label{fig:real-likelihood}
\end{figure}

Figure~\ref{fig:real-likelihood} shows that LM attains lower training negative log-likelihood values within the plotted iteration ranges for both representations. Taken together, the results show that the proposed structured solver can handle the larger Enron tensor and improve held-out predictions in relatively few outer iterations. The higher cost per LM iteration offsets part of this gain, so fewer outer iterations do not imply the shortest runtime. For Enron, AUC helps interpret the results because the large proportion of zeros limits the informativeness of MER alone.

\section{Conclusion}\label{Sec-8}

We developed a structured trust-region LM method for binary tensor completion. Observed-entry contractions and factor-row block preconditioning support matrix-free model solves. Under boundedness of the iterates and sufficient model decrease, the trust-region method converges to stationarity. With exact solves and the original curvature, local quadratic convergence holds under a component-gradient residual error bound for $\delta\in[1,2]$. For the symmetric logistic--LSE family, the analysis gives fourth-order objective-model error and second-order agreement of the reduction ratio. In the HCP and Enron experiments, triple LM gives the best held-out MER and AUC among the tested configurations and uses fewer outer iterations than NCG and L-BFGS.

\bigskip
\noindent\textbf{Data and code availability}~The synthetic example scripts and the reported convergence plots accompany the manuscript. The HCP data are described in \cite{WZD19,WL20}; the Enron experiments use a subset of the public CMU email corpus.

\section*{Acknowledgments}

The work was partly supported by the National Natural Science Foundation of China (12671423,  12671367, 12471351).

\section*{Declarations}
\noindent\textbf{Conflict of interest}~The authors declare no conflicts of interest.

\end{document}